\documentclass[12pt,draftclsnofoot,onecolumn]{IEEEtran}

\usepackage{setspace}
\usepackage{subfig}
\usepackage{epsfig}
\usepackage{cite}
\usepackage{array}
\usepackage{mdwmath}
\usepackage{mdwtab}
\begin{document}

\title{Random Distances Associated with Trapezoids} 

\author{Maryam Ahmadi and Jianping Pan\\
University of Victoria, BC, Canada}

\maketitle

\begin{abstract}
 The distributions of the random distances associated with hexagons, rhombuses and triangles
 have been derived and verified in the existing work. All of these geometric shapes are related 
 to each other and have various applications in wireless
 communications, transportation, etc. Hexagons are widely used to model the cells in cellular
  networks, while trapezoids can 
 be utilized to model the edge users in a cellular network with a hexagonal tessellation.
 In this report, the distributions of the random distances associated with unit trapezoids are derived, 
 when two random points are within a unit trapezoid or in two neighbor unit trapezoids. The 
 mathematical expressions are verified through simulation. Further, we present the polynomial fit 
 for the PDFs, which can be used to simplify the computation.
 \end{abstract}
 

\section{The Problem}

\begin{figure}
  \centering
  \includegraphics[width=0.5\columnwidth]{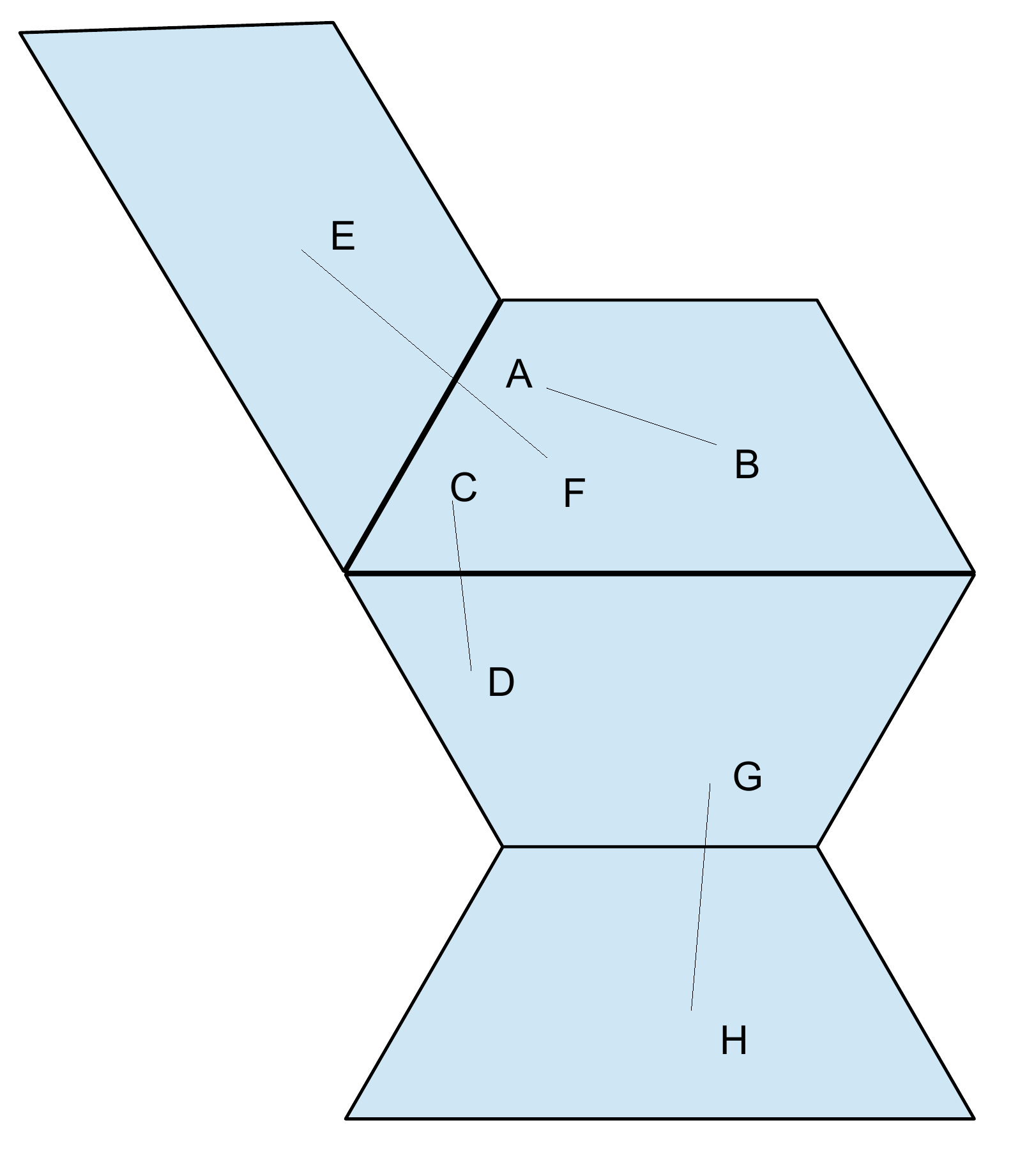}
  \caption{Random Distances Associated with Trapezoids.}
  \label{fig:trapezoids}
\end{figure}

 Denote a ``unit trapezoid" as an isosceles trapezoid where the legs have the same length 
 and equal to the length of the short base, $1$, and the base angles are equal to $\frac{\pi}{3}$.
 The length of the long base is $2$.
 The distributions of the random distances between two points located within a unit trapezoid 
 and in two neighbor unit trapezoids are of our interest. Several cases are considered depending
 on the arrangement of the trapezoids. As illustrated in Fig.~\ref{fig:trapezoids}, points ${\rm A}$ and
 ${\rm B}$ constitute the case where the two random points are located within a unit trapezoid,  
 while ${\rm CD}$, ${\rm EF}$, and ${\rm GH}$ represent the cases where the two random points 
 are located in two neighbor unit trapezoids. For ${\rm CD}$, two neighbor unit trapezoids form 
 a ``unit hexagon" with side length $1$. 
 Two neighbor trapezoids share a leg for ${\rm EF}$. In the last case, ${\rm GH}$, two trapezoids 
 share the short base, forming a concave polygon.   

 We derive the probability density functions (PDFs) and cumulative distribution functions (CDFs) of the $4$ cases in the following section.

\section{Distance Distributions Associated with Unit Trapezoids}\label{sec:result}

\subsection{PDF and CDF of the Random Distances within a Unit Trapezoid (Case $|AB|$)}

A unit trapezoid in fact consists of three adjacent equilateral triangles. Since the distributions of the 
random distances within and between triangles are known~\cite{triangles}, the PDF of the distances 
between two uniformly at random located points within a unit trapezoid can be derived using the  
probabilistic sum, as follows.
 \begin{equation}\label{eq:fd_r_within}
  f_{D_{\rm AB}}(d)=2d\left\{
    \begin{array}{lr}
\frac{8}{9}\left(\frac{\pi}{9\sqrt{3}}+\frac{2}{3}\right)d^2-\frac{80}{27}d+\frac{4\pi}{3\sqrt{3}} 
& 0\leq d\leq \frac{\sqrt{3}}{2}\\

\frac{8}{3}\left(\frac{2}{9\sqrt{3}}d^2+\frac{1}{\sqrt{3}}\right)\sin^{-1}\frac{\sqrt{3}}{2d}+\frac{16}{9}
\left(\frac{1}{3}-\frac{\pi}{9\sqrt{3}}\right)d^2\\
~~~~+\frac{28}{27}\sqrt{4d^2-3}-\frac{80}{27}d & \frac{\sqrt{3}}{2}\leq d\leq 1\\

\frac{16}{9}\left(-\frac{1}{3\sqrt{3}}d^2+\frac{1}{\sqrt{3}}\right)\sin^{-1}\frac{\sqrt{3}}{2d}+
\frac{16\pi}{81\sqrt{3}}d^2-\frac{4}{27}d^2\\
~~~~+\frac{4}{9}\sqrt{4d^2-3}-\frac{32}{27}d+\frac{8\pi}{27\sqrt{3}}-\frac{4}{9} & 1\leq d\leq \sqrt{3} \\

\frac{16}{9\sqrt{3}}\left(\frac{1}{3}d^2+2\right)\sin^{-1}\frac{\sqrt{3}}{d}+\frac{4}{9}\left(\frac{1}{3}-
\frac{4\pi}{9\sqrt{3}}\right)d^2\\
~~~~+\frac{16}{9}\sqrt{d^2-3}-\frac{32}{27}d-\frac{32\pi}{27\sqrt{3}} & \sqrt{3}\leq d\leq 2 \\
      0 & {\rm otherwise}
    \end{array}
  \right..
\end{equation}

The corresponding CDF is 
\begin{equation}\label{eq:Fd_r_within}
  F_{D_{\rm AB}}(d)=\left\{
    \begin{array}{lr}
0&d<0\\
\frac{4}{9}\left(\frac{2}{3}+\frac{\pi}{9\sqrt{3}}\right)d^4-\frac{160}{81}d^3+\frac{
4\pi}{3\sqrt{3}}d^2 & 0\leq d\leq \frac{\sqrt{3}}{2}\\

\frac{8}{27\sqrt{3}}\left(d^2+9\right)d^2\sin^{-1}\frac{\sqrt{3}}{2d}+\frac{8}{9}\left(\frac{1}{3}-
\frac{\pi}{9\sqrt{3}}\right)d^4\\
~~~~-\frac{160}{81}d^3+\frac{58d^2+15}{81}\sqrt{4d^2-3} & \frac{\sqrt{3}}{2}\leq d\leq 1\\

-\frac{8}{27\sqrt{3}}(d^2-6)d^2\sin^{-1}\frac{\sqrt{3}}{2d}+\left(\frac{8\pi}{81\sqrt{3}}-
\frac{2}{27}\right)d^4\\
~~~~-\frac{64}{81}d^3+\left(\frac{8\pi}{27\sqrt{3}}-\frac{4}{9}\right)d^2\\
~~~~+\frac{22d^2+15}{81}\sqrt{4d^2-3}+0.0735 & 1\leq d\leq \sqrt{3} \\

\frac{8}{27\sqrt{3}}(d^2+12)d^2\sin^{-1}\frac{\sqrt{3}}{d}+\frac{2}{9}\left(\frac{1}{3}-
\frac{4\pi}{9\sqrt{3}}\right)d^4\\
~~~~-\frac{64}{81}d^3-\frac{32\pi}{27\sqrt{3}}d^2+\frac{104d^2+48}{81}\sqrt{d^2-3}\\
~~~~+0.4074 & \sqrt{3}\leq d\leq 2 \\
1 & d>2
    \end{array}
 \right..
\end{equation}

\subsection{PDF and CDF of the Distances between Two Uniformly at Random Located Points
 in Neighbor Unit Trapezoids (Case $|CD|$)}

In this case, the two neighbor trapezoids form a hexagon. Since the distributions of the random 
distances within a unit hexagon have been derived and verified in~\cite{hexagons}, the PDF  
of the distances between two uniformly at random located nodes in two neighbor unit trapezoids 
can be derived as
 \begin{equation}\label{eq:fd_r_hex}
  f_{D_{\rm CD}}(d)=2d\left\{
    \begin{array}{lr}
    
-\frac{4}{27}\left(1+\frac{5\pi}{3\sqrt{3}}\right)d^2+\frac{32}{27}d & 0\leq d\leq
\frac{\sqrt{3}}{2}\\

-\frac{1}{3\sqrt{3}}\left(\frac{16}{9}d^2+8\right)\sin^{-1}\frac{\sqrt{3}}{2d}+\frac{4}{27}\left(\frac{\pi}
{3\sqrt{3}}-1\right)d^2\\
~~~~-\frac{28}{27}\sqrt{4d^2-3}+\frac{32}{27}d+\frac{4\pi}{3\sqrt{3}} & \frac{\sqrt{3}}{2}\leq d\leq 1\\

 -\frac{4}{9\sqrt{3}}\left(\frac{8}{3}d^2+10\right)\sin^{-1}\frac{\sqrt{3}}{2d}+\frac{4}{27}\left(\frac{5\pi}
 {3\sqrt{3}}+1\right)d^2\\
 ~~~~-\frac{16}{9}\sqrt{4d^2-3}+\frac{32}{27}d+\frac{52\pi}{27\sqrt{3}}+\frac{4}{9}& 1\leq d\leq \sqrt{3} \\

\frac{8}{9\sqrt{3}}\left(\frac{1}{3}d^2+8\right)\sin^{-1}\frac{\sqrt{3}}{d}-\frac{8}{27}\left(\frac{\pi}
{3\sqrt{3}}+2\right)d^2\\
~~~~+\frac{8}{3}\sqrt{d^2-3}+\frac{32}{27}d-\frac{64\pi}{27\sqrt{3}}-\frac{8}{3} & \sqrt{3}\leq d\leq 2 \\
      0 & {\rm otherwise}
    \end{array}
 \right..
\end{equation}

The corresponding CDF is 
\begin{equation}\label{eq:Fd_r_para}
F_{D_{\rm CD}}(d)=\left\{
    \begin{array}{lr}    
    0 & d<0 \\
-\frac{2}{27}\left(\frac{5\pi}{3\sqrt{3}}+1\right)d^4+\frac{64}{81}d^3 & 0\leq d\leq \frac{\sqrt{3}}{2}\\

-\frac{8}{27\sqrt{3}}\left(d^2+9\right)d^2\sin^{-1}\frac{\sqrt{3}}{2d}+\frac{2}{27}\left(\frac{\pi}
{3\sqrt{3}}-1\right)d^4\\
~~~~+\frac{64}{81}d^3+\frac{4\pi}{3\sqrt{3}}d^2-\frac{58d^2+15}{81}\sqrt{4d^2-3} & \frac{\sqrt{3}}{2}\leq d\leq 1\\

-\frac{8}{27\sqrt{3}}(2d^2+15)d^2\sin^{-1}\frac{\sqrt{3}}{2d}+\frac{2}{27}\left(\frac{5\pi}{3\sqrt{3}}+
1\right)d^4\\
~~~~+\frac{64}{81}d^3+\left(\frac{52\pi}{27\sqrt{3}}+\frac{4}{9}\right)d^2-\frac{100d^2+24}{81}
\sqrt{4d^2-3}\\
~~~~+0.0359 & 1\leq d\leq \sqrt{3} \\

\frac{4}{27\sqrt{3}}(d^2+48)d^2\sin^{-1}\frac{\sqrt{3}}{d}-\frac{4}{27}\left(\frac{\pi}{3\sqrt{3}}+2\right)d^4\\
~~~~+\frac{64}{81}d^3-\left(\frac{64\pi}{27\sqrt{3}}+\frac{8}{3}\right)d^2+\frac{148d^2+168}{81}\sqrt{d^2-3}\\
~~~~+0.7026 & \sqrt{3}\leq d\leq 2 \\
1 & d>2
    \end{array}
 \right..
\end{equation}

\subsection{PDF and CDF of the Distances between Two Uniformly at Random Located Points in Neighbor 
Unit Trapezoids (Case $|EF|$)}

Knowing the PDF of the random distances between rhombuses~\cite{rhombuses} and within a unit 
trapezoid, the PDF of the random distances between two uniformly distributed points, one in each 
of the neighbor unit trapezoids in this case, is

 \begin{equation}\label{eq:fd_r_diag1}
  f_{D_{\rm EF}}(d)=2d\left\{
    \begin{array}{lr}
-\frac{4}{27}\left(\frac{2\pi}{3\sqrt{3}}+1\right)d^2+\frac{16}{27}d & 0\leq d\leq \frac{\sqrt{3}}{2}\\

-\frac{8}{9\sqrt{3}}\left(\frac{2}{3}d^2+1\right)\sin^{-1}\frac{\sqrt{3}}{2d}+\frac{4}{27}\left(\frac{4\pi}{3\sqrt{3}}
-1\right)d^2\\
~~~~+\frac{16}{27}d+\frac{4\pi}{9\sqrt{3}}-\frac{4}{9}\sqrt{4d^2-3} & \frac{\sqrt{3}}{2}\leq d\leq 1\\

\frac{8}{9\sqrt{3}}\left(\frac{1}{3}d^2+1\right)\sin^{-1}\frac{\sqrt{3}}{2d}-\frac{4}{27}\left(\frac{2\pi}{3\sqrt{3}}
+1\right)d^2\\
~~~~+\frac{10}{27}\sqrt{4d^2-3}-\frac{4\pi}{27\sqrt{3}}-\frac{2}{9} & 1\leq d\leq \sqrt{3}\\

 \frac{8}{9\sqrt{3}}\left(\frac{1}{3}d^2-2\right)\sin^{-1}\frac{\sqrt{3}}{2d}
  -\frac{8}{9\sqrt{3}}\left(\frac{1}{3}d^2+2\right)\sin^{-1}\frac{\sqrt{3}}{d}\\
 ~~~~+\frac{4}{27}\left(\frac{\pi}{3\sqrt{3}}+2\right)d^2
-\frac{8}{9}\sqrt{d^2-3}-\frac{14}{27}\sqrt{4d^2-3}\\
 ~~~~+\frac{32\pi}{27\sqrt{3}}+\frac{10}{9}& \sqrt{3}\leq d\leq 2\\

 \frac{8}{9\sqrt{3}}\left(\frac{1}{3}d^2-2\right)\sin^{-1}\frac{\sqrt{3}}{2d}
  +\frac{16}{9\sqrt{3}}\sin^{-1}\frac{\sqrt{3}}{d}
 -\frac{4}{27}\left(\frac{\pi}{3\sqrt{3}}-1\right)d^2\\
~~~~-\frac{14}{27}\sqrt{4d^2-3}+\frac{16}{27}\sqrt{d^2-3}
 +\frac{2}{9}
 &2\leq d\leq \sqrt{7} \\

 \frac{8}{9\sqrt{3}}\left(-\frac{1}{3}d^2+4\right)\sin^{-1}\frac{\sqrt{3}}{d}
 +\frac{4}{27}\left(\frac{\pi}{3\sqrt{3}}-1\right)d^2\\
~~~~+\frac{8}{9}\sqrt{d^2-3}
 -\frac{16\pi}{27\sqrt{3}} -\frac{8}{9}
 & \sqrt{7} \leq d \leq 2\sqrt{3} \\
      0 & {\rm otherwise}
    \end{array}
  \right..
\end{equation}

The corresponding CDF is 
\begin{equation}\label{eq:Fd_r_diag1}
  F_{D_{\rm EF}}(d)=\left\{
    \begin{array}{lr}
    0 & d<0\\
-\frac{2}{27}\left(\frac{2\pi}{3\sqrt{3}}+1\right)d^4+\frac{32}{81}d^3 & 0\leq d\leq \frac{\sqrt{3}}{2}\\

-\frac{8}{27\sqrt{3}}\left(d^2+3\right)d^2\sin^{-1}\frac{\sqrt{3}}{2d}+\frac{2}{27}\left(\frac{4\pi}{3\sqrt{3}}
-1\right)d^4\\
~~~~+\frac{32}{81}d^3+\frac{4\pi}{9\sqrt{3}}d^2-\frac{26d^2+3}{81}\sqrt{4d^2-3} & \frac{\sqrt{3}}{2}\leq d\leq 1\\

\frac{4}{27\sqrt{3}}(d^2+6)d^2\sin^{-1}\frac{\sqrt{3}}{2d}-\frac{2}{27}\left(\frac{2\pi}{3\sqrt{3}}+1\right)d^4\\
~~~~-\left(\frac{4\pi}{27\sqrt{3}}+\frac{2}{9}\right)d^2+\frac{42d^2+9}{162}\sqrt{4d^2-3}-0.0561
& 1\leq d\leq \sqrt{3} \\

\frac{4}{27\sqrt{3}}d^2\left((d^2-12)\sin^{-1}\frac{\sqrt{3}}{2d}-(d^2+12)\sin^{-1}\frac{\sqrt{3}}{d}\right)\\
~~~~+\frac{2}{27}\left(\frac{\pi}{3\sqrt{3}}+2\right)d^4+\frac{2}{9}\left(\frac{16\pi}{3\sqrt{3}}+5\right)d^2\\
~~~~-\frac{54d^2+27}{162}\sqrt{4d^2-3}-\frac{52d^2+24}{81}\sqrt{d^2-3}-0.0561 & \sqrt{3}\leq d\leq 2 \\

\frac{4}{27\sqrt{3}}(d^2-12)d^2\sin^{-1}\frac{\sqrt{3}}{2d}+\frac{16}{9\sqrt{3}}d^2\sin^{-1}\frac{\sqrt{3}}{d}\\
~~~~-\frac{2}{27}\left(\frac{\pi}{3\sqrt{3}}-1\right)d^4+\frac{2}{9}d^2-\frac{54d^2+27}{162}\sqrt{4d^2-3}\\
~~~~+\frac{32d^2+48}{81}\sqrt{d^2-3}-0.3528 & 2\leq d\leq \sqrt{7} \\ 

-\frac{4}{27\sqrt{3}}(d^2-24)d^2\sin^{-1}\frac{\sqrt{3}}{d}
+\frac{2}{27}\left(\frac{\pi}{3\sqrt{3}}-1\right)d^4\\
~~~~-\frac{8}{9}\left(\frac{2\pi}{3\sqrt{3}}+1\right)d^2+\frac{44d^2+120}{81}\sqrt{d^2-3}-1.6677 & 
\sqrt{7}\leq d\leq 2\sqrt{3} \\
1&d>2\sqrt{3}
    \end{array}
  \right..
\end{equation}

\subsection{PDF and CDF of the Distances between Two Uniformly at Random Located Points in Neighbor 
Unit Trapezoids (Case $|GH|$)}

The PDF of the random distances between two nodes uniformly deployed in two neighbor unit
 trapezoids (case GH), is
 \begin{equation}\label{eq:fd_r_diag2}
  f_{D_{\rm GH}}(d)=2d\left\{
    \begin{array}{lr}
\frac{4}{27}\left(\frac{\pi}{3\sqrt{3}}-1\right)d^2+\frac{16}{27}d
 & 0\leq d\leq\frac{\sqrt{3}}{2}\\

 \frac{8}{9\sqrt{3}}\left(\frac{2}{3}d^2-1\right)\sin^{-1}\frac{\sqrt{3}}{2d}
  -\frac{4}{27}\left(\frac{5\pi}{3\sqrt{3}}+1\right)d^2\\
~~~~-\frac{4}{27}\sqrt{4d^2-3}
 +\frac{16}{27}d+\frac{4\pi}{9\sqrt{3}}
 & \frac{\sqrt{3}}{2} \leq d \leq1\\

 \frac{8}{3\sqrt{3}}\left(\frac{1}{9}d^2-1\right)\sin^{-1}\frac{\sqrt{3}}{2d}
  -\frac{4}{9}\left(\frac{\pi}{3\sqrt{3}}-1\right)d^2\\
 ~~~~-\frac{22}{27}\sqrt{4d^2-3}
 +\frac{28\pi}{27\sqrt{3}}+\frac{2}{3}
 & 1\leq d\leq \sqrt{3} \\

 -\frac{8}{9\sqrt{3}}\left(\frac{1}{3}d^2-2\right)\sin^{-1}\frac{\sqrt{3}}{2d}
  -\frac{8}{9\sqrt{3}}\left(\frac{1}{3}d^2-2\right)\sin^{-1}\frac{\sqrt{3}}{d}\\
 ~~~~+\frac{8}{27}\left(\frac{\pi}{3\sqrt{3}}-1\right)d^2
+\frac{14}{27}\sqrt{4d^2-3}+\frac{8}{27}\sqrt{d^2-3}\\
 ~~~~-\frac{16\pi}{27\sqrt{3}}-\frac{10}{9}
 & \sqrt{3}\leq d\leq \sqrt{7}  \\
      0 & {\rm otherwise}
    \end{array}
  \right..
\end{equation}

The corresponding CDF is 
 \begin{equation}\label{eq:Fd_r_diag2}
  F_{D_{\rm GH}}(d)=\left\{
    \begin{array}{lr}
    0&d<0\\
\frac{2}{27}\left(\frac{\pi}{3\sqrt{3}}-1\right)d^4+\frac{32}{81}d^3 & 0\leq d\leq \frac{\sqrt{3}}{2}\\

\frac{8}{27\sqrt{3}}\left(d^2-3\right)d^2\sin^{-1}\frac{\sqrt{3}}{2d}-\frac{2}{27}\left(\frac{5\pi}{3\sqrt{3}}
+1\right)d^4\\
~~~~+\frac{32}{81}d^3+\frac{4\pi}{9\sqrt{3}}d^2-\frac{2d^2+3}{27}\sqrt{4d^2-3} & \frac{\sqrt{3}}{2}\leq d\leq 1\\

\frac{4}{27\sqrt{3}}(d^2-18)d^2\sin^{-1}\frac{\sqrt{3}}{2d}-\frac{2}{9}\left(\frac{\pi}{3\sqrt{3}}-1\right)d^4\\
~~~~+\frac{2}{3}\left(\frac{14\pi}{9\sqrt{3}}+1\right)d^2-\frac{86d^2+39}{162}\sqrt{4d^2-3}+0.0176
& 1\leq d\leq \sqrt{3} \\

-\frac{4}{27\sqrt{3}}(d^2-12)d^2\left(\sin^{-1}\frac{\sqrt{3}}{2d}+\sin^{-1}\frac{\sqrt{3}}{d}\right)+
\frac{4}{27}\left(\frac{\pi}{3\sqrt{3}}-1\right)d^4\\
~~~~-\frac{2}{9}\left(\frac{8\pi}{3\sqrt{3}}+5\right)d^2+\frac{54d^2+27}{162}\sqrt{4d^2-3}+
\frac{12d^2+72}{81}\sqrt{d^2-3}\\
~~~~-0.3157 & \sqrt{3}\leq d\leq \sqrt{7} \\
1&d>\sqrt{7}
    \end{array}
  \right..
\end{equation}

The results presented above are for the ``unit" trapezoids. However, the results can be easily extended for scaled 
trapezoids. Assume each side of the unit trapezoid is scaled by a non-negative parameter $s$, then,
\begin{equation}
 F_{sD}(d)=P(sD\leq d)=P(D\leq \frac{d}{s})=F_D(\frac{d}{s}). \nonumber
\end{equation}
Therefore,
\begin{equation}\label{eq:scale}
 f_{sD}(d)=F'_D(\frac{d}{s})=\frac{1}{s}f_D(\frac{d}{s}).
\end{equation}

\section{Verification}

\begin{figure}
  \centering
  \includegraphics[width=0.55\columnwidth]{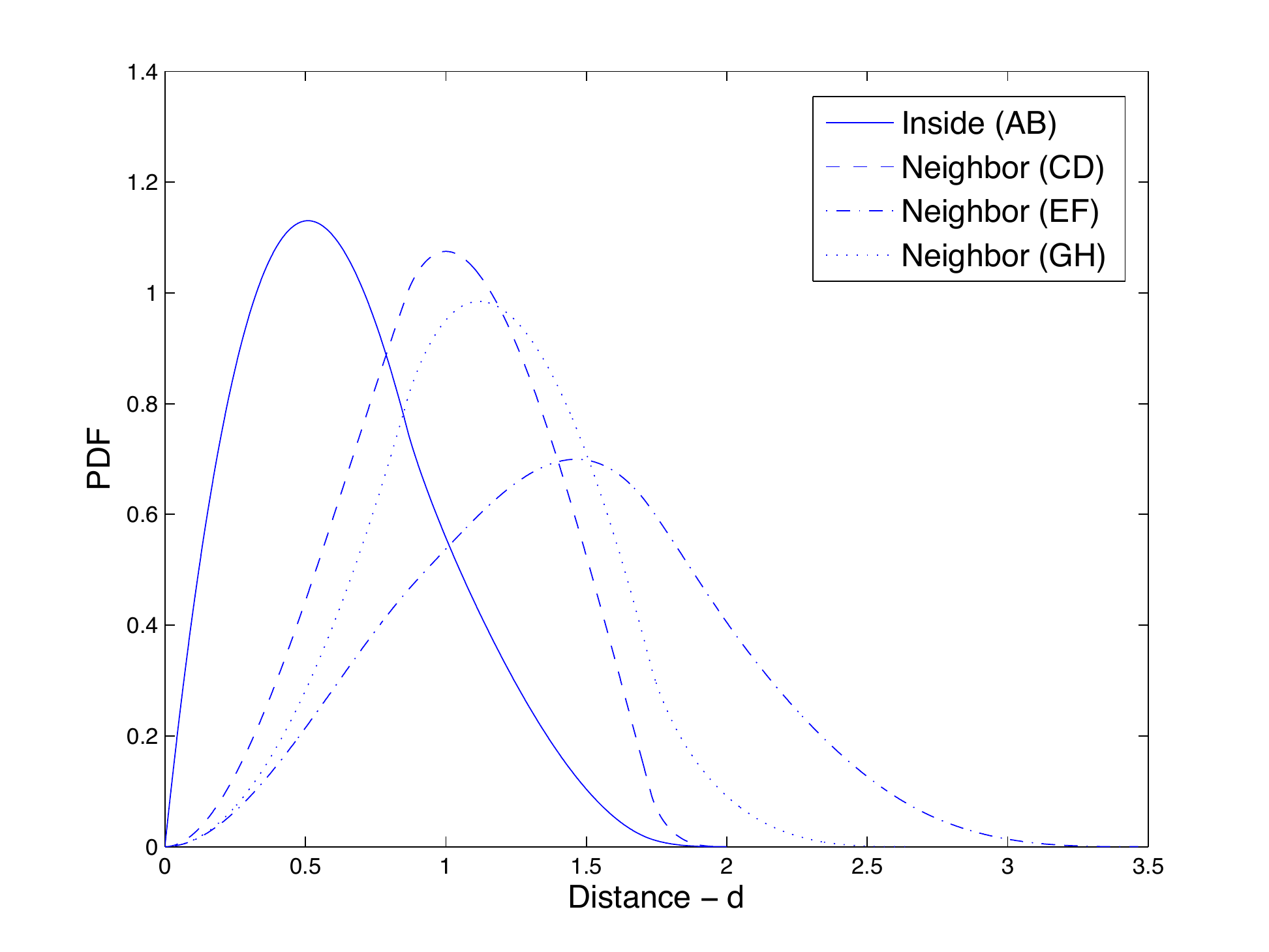}
  \caption{Distributions of the Random Distances Associated with Unit Trapezoids.}
  \label{fig:pdfs}
\end{figure}

\begin{figure}
  \centering
  \includegraphics[width=0.55\columnwidth]{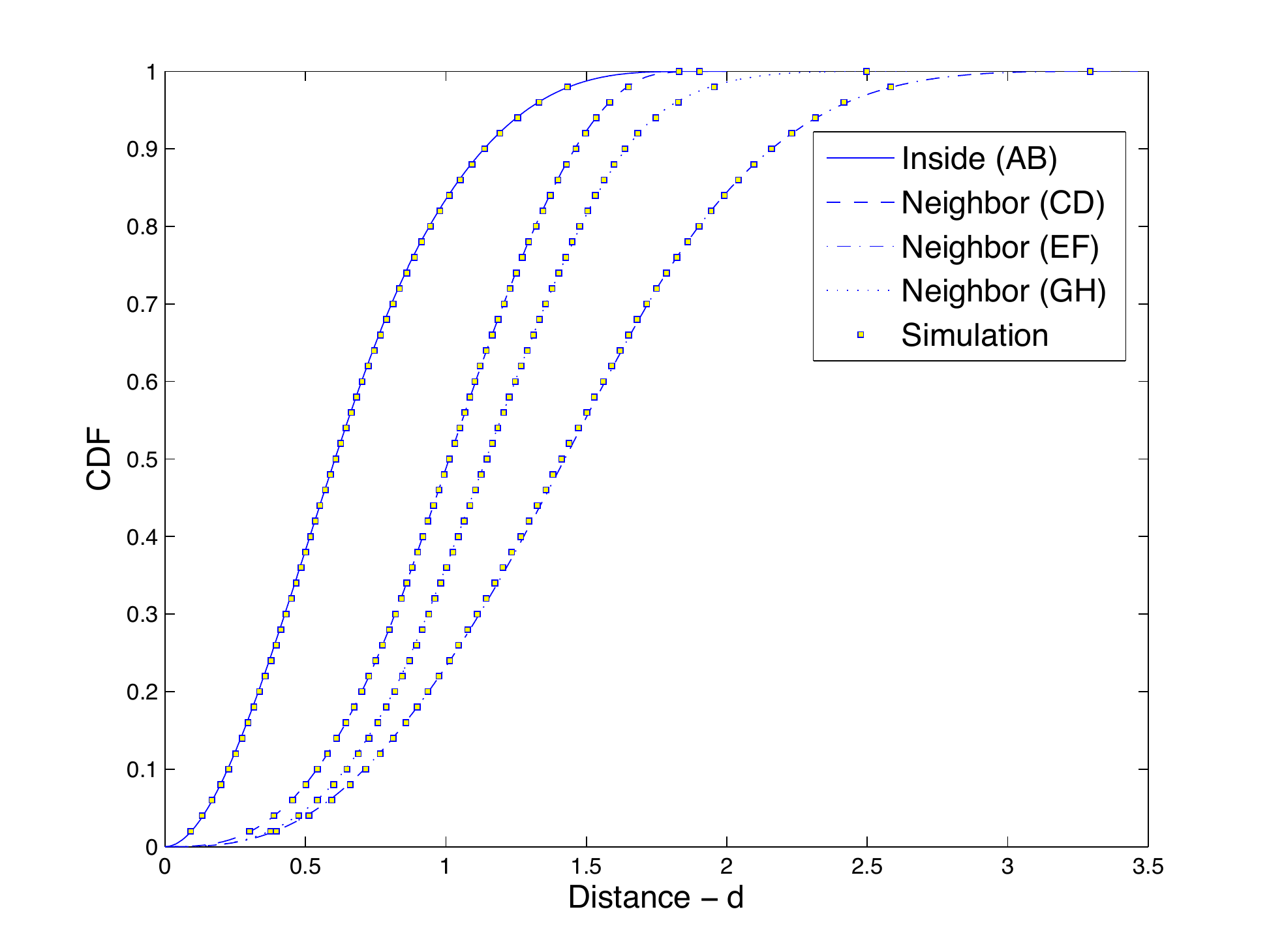}
  \caption{Comparing Simulation and Analytical results.}
  \label{fig:cdfs}
\end{figure}

Figure~\ref{fig:pdfs} demonstrates the PDF of the random distances associated with trapezoids, 
including the random distances within a unit trapezoid and between neighbor unit trapezoids.

Figure~\ref{fig:cdfs} shows the comparison between the simulation and analytical results of 
the CDFs. The simulation is done in Matlab, generating $10,000$ pairs of random points.
As the figure suggests, there is a close match between the mathematical and simulation results, 
which verifies the accuracy of our analytical results.

\section{Polynomial Fits}

\begin{table}
  \caption{Coefficients of the Polynomial Fit and the Norm of Residuals (NormR)}
  \centering
  \begin{tabular}{|c||c|c|}
    \hline
    PDF & Polynomial Coefficients & NormR \\ \hline \hline
    & $\left[-6.3419~~72.4436~~-356.68~~988.25~~-1686.9\right.$ & \\  
    $f_{D_{\rm AB}}(d)$ &
$1824.8~~-1239.4~~504.30~~-109.17~~9.6850~~-5.1833$ & $0.1537$ \\ 
     & $\left.4.7245~~0.0021\right]$ & 
     \\\hline
     & $
\left[10.55774~~-124.22~~632.96~~-1829.6~~3302.5\right.$ & \\
     $f_{D_{\rm CD}}(d)$ &
$-3864.6~~2953.0~~-1452.6~~443.86~~-80.1282~~9.6421$& $0.1692$ \\
     & $\left.-0.2762 0.0024\right]$ &
\\\hline
& $
\left[0.0196~~-0.3853~~3.2716~~-15.6480~~46.2235\right.$ & \\ 
     $f_{D_{\rm EF}}(d)$ &
$-87.0513 104.33e -77.1478 32.4560 -6.2059$ & $0.2405$ \\
     & $\left.0.6906 -0.0154\right]$ & \\
\hline
   & $
\left[0.6286~~-16.1105~~1.9174\times 10^3~~-1.4062\times 10^3\right.$ & \\
     $f_{D_{\rm GH}}(d)$ &
$7.1067\times 10^3~~-2.6205\times 10^4~~7.2760\times 10^4~~-1.5469\times 10^5$ & $0.1622$ \\
     & $2.5338\times 10^5~~-3.1899\times 10^5~~3.0528\times 10^5~~-2.1736e\times 10^5$ & \\ 
     & $1.1081\times 10^5~~-3.7549\times 10^4~~6.9028\times 10^3~~41.6579$ & \\
     & $\left.-3.2855\times 10^2~~69.7237~~-4.9233~~0.2052~~-0.0016\right]$ & \\
\hline
  \end{tabular}
  \label{tab:poly}
\end{table}

\begin{figure}
\centering
  \subfloat[Within a Unit Trapezoid (AB)]{\includegraphics[width=0.5\columnwidth]{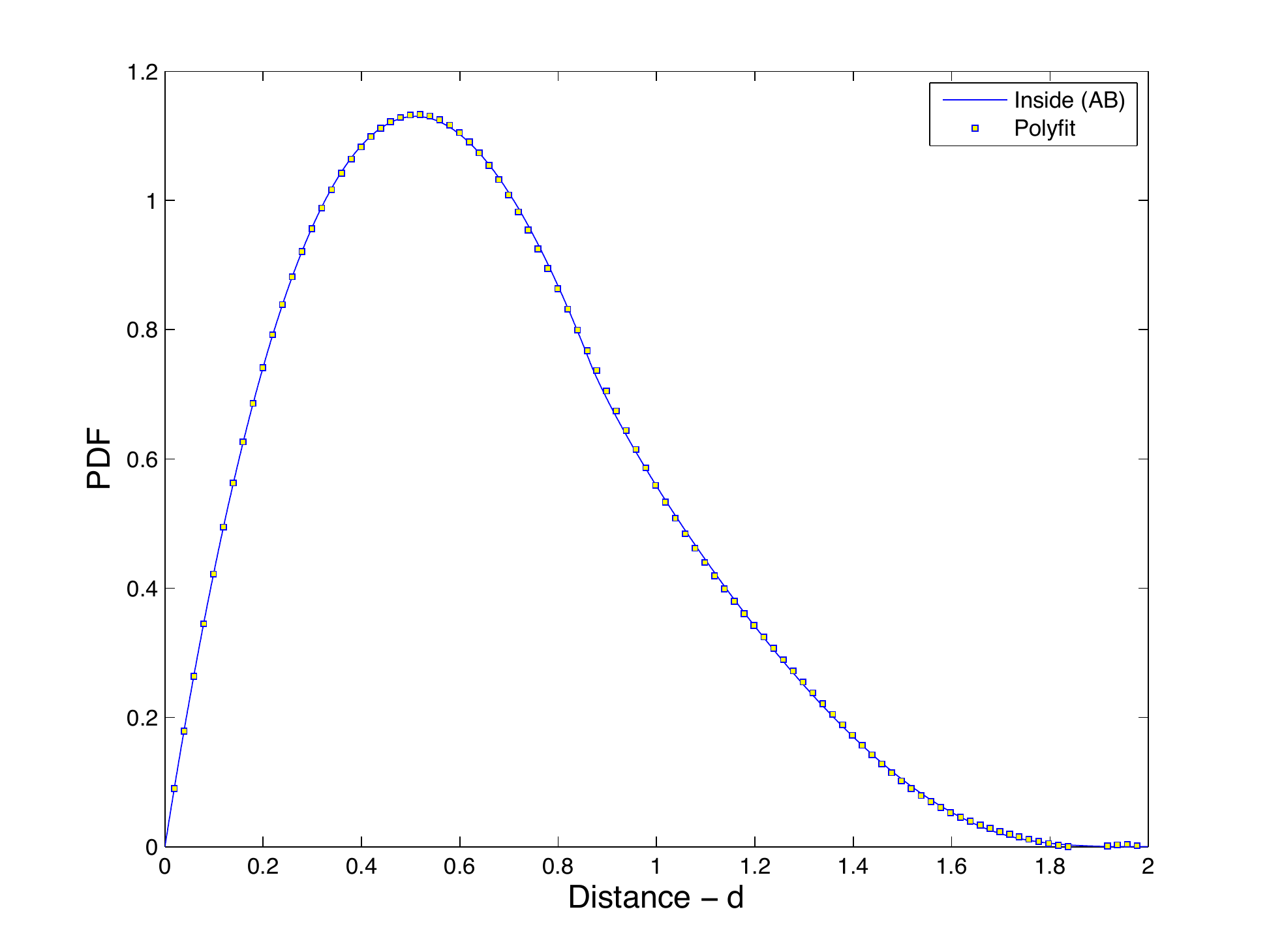}}
  \subfloat[Between Two Neighbor Unit Trapezoids (CD)]{\includegraphics[width=0.5\columnwidth]{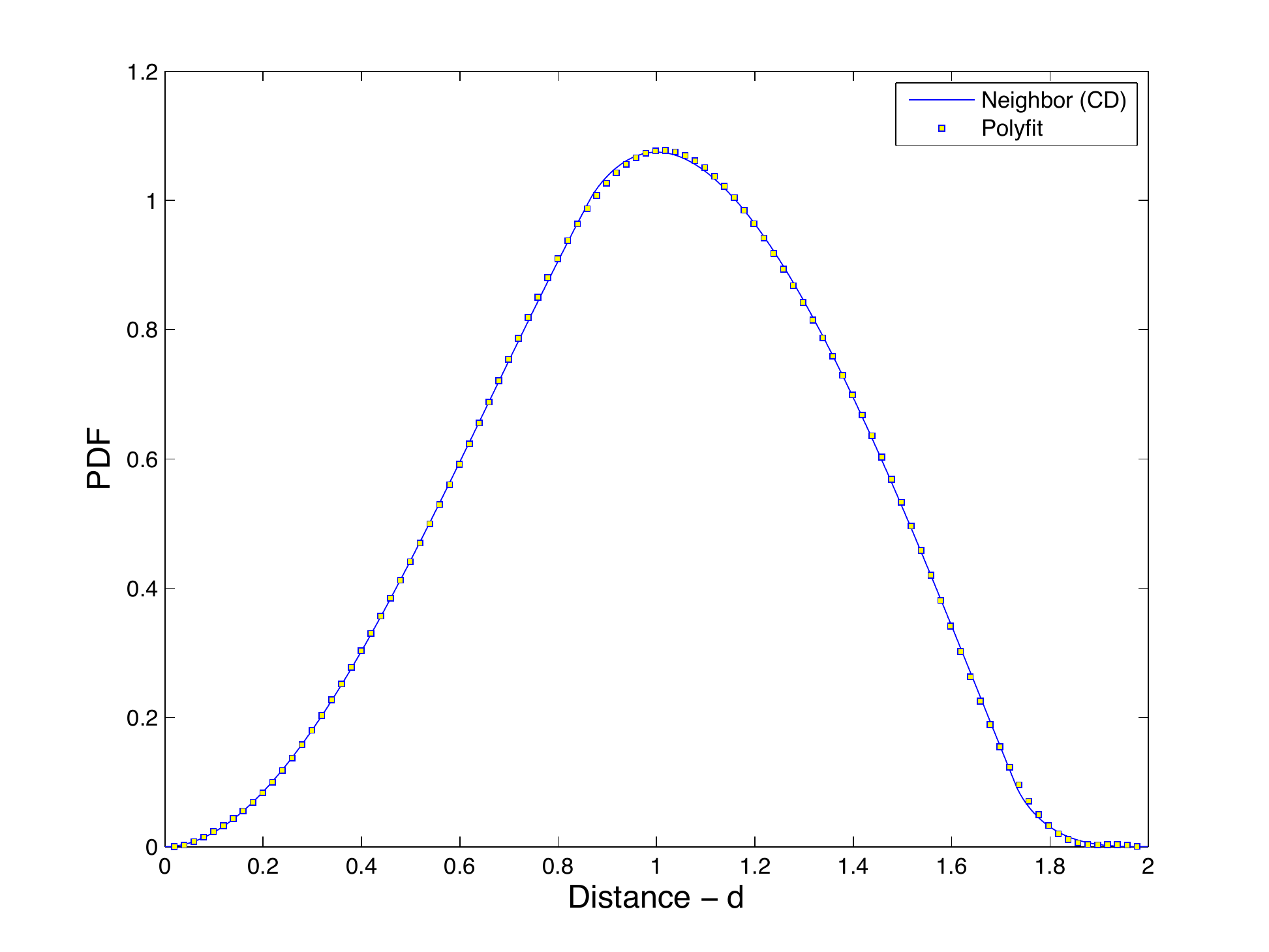}}\\
  \subfloat[Between Two Neighbor Unit Trapezoids (EF)]{\includegraphics[width=0.5\columnwidth]{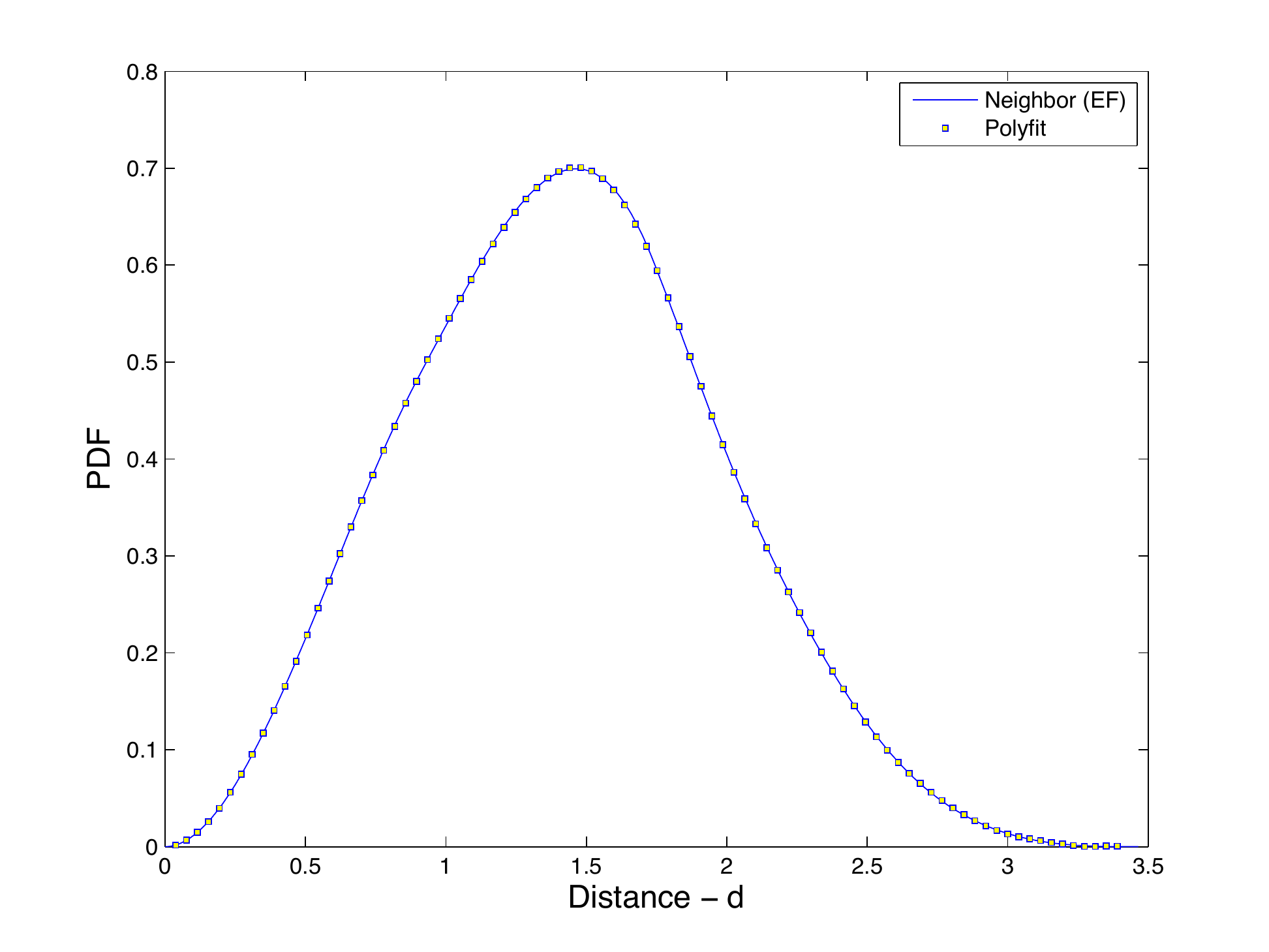}}
   \subfloat[Between Two Neighbor Unit Trapezoids (GH)]{\includegraphics[width=0.5\columnwidth]{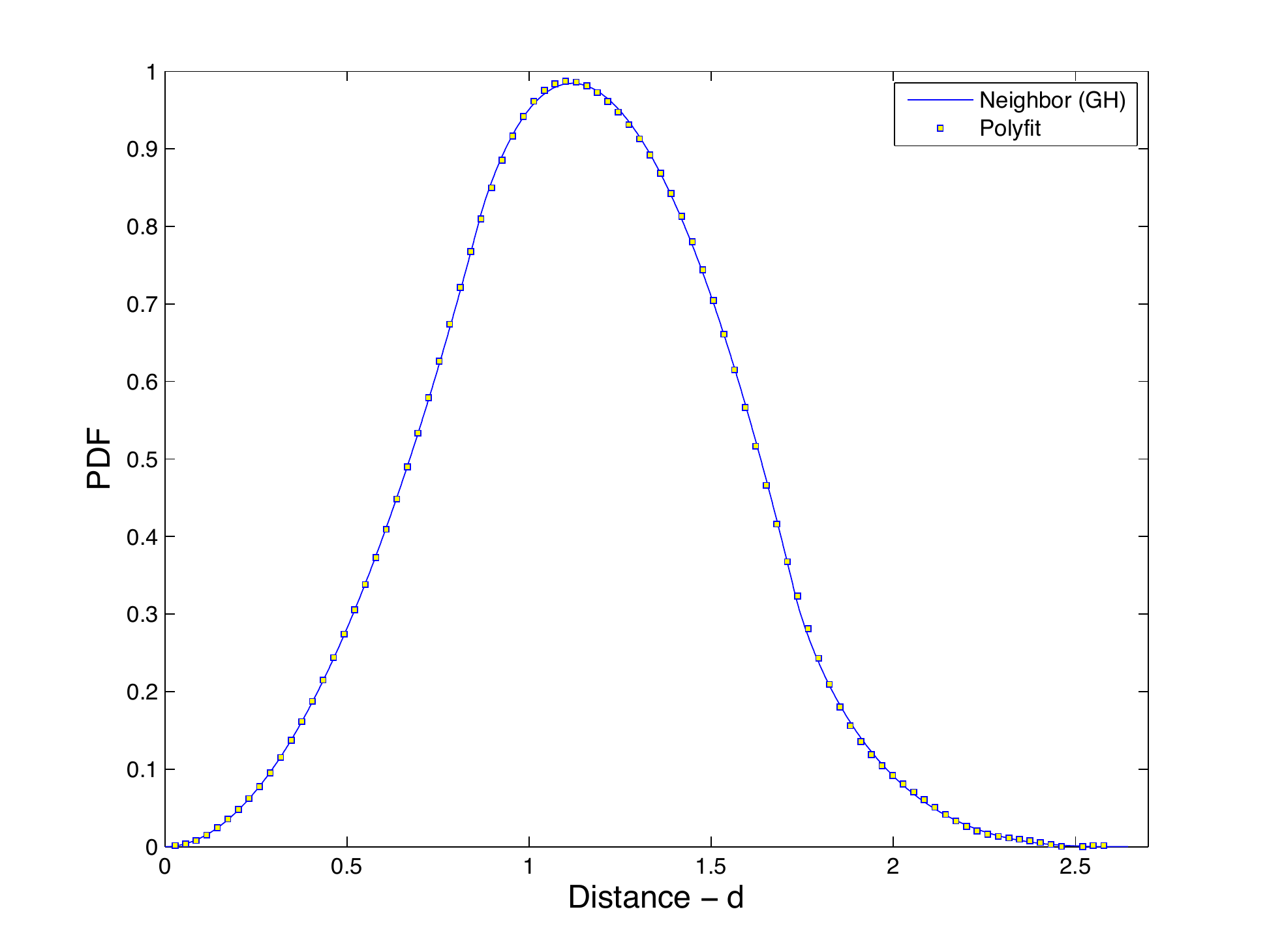}}
  \caption{Polynomial Fits.}
  \label{fig:trapezoid_poly}
\end{figure}

Figure~\ref{fig:trapezoid_poly} shows the PDFs of the random distances within a unit trapezoid as well 
as between two neighbor unit trapezoids, along with the fitted polynomial functions. The polynomial  
coefficients for the PDFs of the random distances associated with trapezoids, along with their norm  
of residual, are presented in Table~\ref{tab:poly}.
In Fig.~\ref{fig:trapezoid_poly}, the approximated polynomial functions are shown with dots over the 
original PDF curves. One can observe that the polynomial functions demonstrate a good match with 
the original PDFs. The polynomial functions are accurate and can be used instead of the original  
PDFs to reduce the computational complexity.

\section{Conclusions}
\label{sec:conclude}

In this report, the closed-form expressions for the random distances within a unit trapezoid 
and between two neighbor unit trapezoids are given. 
The analytical results are verified through simulation. In addition, the polynomial fits for the PDFs 
are obtained. These polynomials offer a good fit and can be used instead of the original PDFs 
to reduce the computational complexity.

\section*{Acknowledgment}
This work is supported in part by the NSERC, CFI and BCKDF. The authors would like to thank
 Lei Zhang, Tianming Wei and Fei Tong for their help.

\end{document}